\documentclass[11pt]{article}

\usepackage{tgpagella}
\usepackage[utf8]{inputenc}
\usepackage[T1]{fontenc}
\usepackage{amssymb}
\usepackage{amsfonts}
\usepackage{amsmath}
\usepackage{amsthm}
\usepackage{mathrsfs}
\usepackage{comment}
\usepackage{xcolor}
\usepackage{hyperref}

\newtheorem{theorem}{Theorem}
\newtheorem{proposition}[theorem]{Proposition}
\newtheorem{lemma}[theorem]{Lemma}

\newtheorem{corollary}[theorem]{Corollary}

\theoremstyle{definition}
\newtheorem{definition}[theorem]{Definition}

\newtheorem{example}[theorem]{Example}

\newtheorem{question}[theorem]{Question}

\newcommand{\PA}{\textnormal{PA}}
\newcommand{\FT}{\textnormal{FT}}

\newcommand{\set}[2]{\lbrace #1 \ \mid \ #2 \rbrace }
\newcommand{\res}{\upharpoonright}
\newcommand{\CT}{\textnormal{CT}}
\newcommand{\df}[1]{\textbf{#1}}
\newcommand{\num}[1]{\underline{#1}}

\newcommand{\INT}{\textnormal{INT}}
\newcommand{\INTColl}{\textnormal{INTColl}}

\newcommand{\Coll}{\textnormal{Coll}}
\newcommand{\ElDiag}{\textnormal{ElDiag}}
\newcommand{\val}[1]{{#1}^{\circ}}

\newcommand{\tuple}[1]{\langle #1 \rangle}

\newcommand{\dpt}{\textnormal{dp}}

\newcommand{\dom}{\textnormal{dom}}

\newcommand{\id}{\textnormal{id}}

\newcommand{\LPA}{\mathscr{L}_{\PA}}
\newcommand{\form}{\textnormal{Form}}

\newcommand{\Term}{\textnormal{Term}}
\newcommand{\Sent}{\textnormal{Sent}}

\newcommand{\ClTerm}{\textnormal{ClTerm}}

\newcommand{\ClTermSeq}{\textnormal{ClTermSeq}}

\newcommand{\Asn}{\textnormal{Asn}}
\newcommand{\Var}{\textnormal{Var}}

\title{Truth and collection} 
\author{Bartosz Wcisło\footnote{Institute of Philosophy, University of Gdańsk.}}

\begin{document}
	
	\maketitle

\begin{abstract}
	Answering a question of Kaye, we show that the compositional truth theory with a full collection scheme is conservative over Peano Arithmetic. We demonstrate it by showing that countable models of compositional truth which satisfy the internal induction or collection axioms can be end-extended to models of the respective theory.
\end{abstract}	

\section{Introduction}

The area of axiomatic truth theories, studies extensions of foundational axiomatic systems such as Peano arithmetic ($\PA$), elementary arithmetic, or Zermelo--Fr\"ankel set theory with axioms expressing that a fresh unary predicate $T(x)$ defines the set of true sentences. 

One of the most basic such theories is the theory of compositional truth over $\PA$, called $\CT^-$.\footnote{In the literature, this theory is also called $\CT \res$, $\CT[\PA]$, $\PA^{\FT}$. There are also variants of these names for the related concept of satisfaction.} Its axioms state that the truth predicate satisfies Tarski's compositional conditions for the arithmetical language. It turns out that the Tarski's axioms themselves constitute a conservative extension of the base theory, even though usually a theory which can formulate a truth predicate for its subtheory is significantly stronger.

In a line of research started by \cite{cies} and \cite{cieslinskict0} and discussed thoroughly in \cite{cies_ksiazka} the following question has been explored: What axioms can be added to $\CT^-$ so as to make the resulting theory nonconservative? The ``line'' dividing conservative truth principles from the non-conservative ones ha been dubbed by Ali Enayat \textit{the Tarski Boundary} and the problem of systematically investigating the conservativity of truth-theoretic extensions has been called \textit{the Tarski Boundary Problem}.

The principles under investigation come roughly in two flavours: either they are of purely truth-theoretic nature, for instance stating that a disjunction of an arbitrary finite size is true iff one of the disjunct is, or they are general principles studied in the context of arithmetical theories, for instance, fragments of the induction scheme. 

It turned out that the investigated principles turn out to be either conservative or equivalent exactly to the scheme of $\Delta_0$-induction for the full language including the truth predicate. Thus, the Tarski Boundary Problem seems to admit a much more structured answer than one could expect. This is a very surprising phenomenon, since \textit{prima facie} the principles under investigation seem completely unrelated.

The topic of axiomatic truth has been classically studied as a part of the theory of models of $\PA$, as a study of so-called \emph{satisfaction classes}. A satisfaction class $S$ in a model $M$ of $\PA$ is a subset of $M^2$ which satisfies Tarski's compositional conditions. They were introduced in the PhD thesis of \cite{krajewski_phd}, preceded by related ideas in \cite{robinson_satisfaction}. They are supposed to explore how  we can endow nonstandard elements of a model $M$ which are viewed from its point of view as arithmetical formulae with well-defined semantics. From this point of view more important than conservativity of truth-theoretic principles is the question whether in a given model of $\PA$ one can find a satisfaction class enjoying additional good properties. 

One of the most pressing questions concerning the Tarski Boundary problem which was left open has been posed by Kaye.\footnote{It was stated explicitly in the presentation \cite{kaye-slides}, but we are not aware of the author stating it in print.
} 
It asks, whether in any countable recursively saturated model of $\PA$, we can find  a full satisfaction class satisfying collection. Equivalently: is $\CT^-$ with the full collection scheme for the extended language, but with no induction whatsoever, a conservative extension of $\PA$? Results about the purely arithmetical counterpart of this question supported the intuition that the answer should be positive: If we add the full collection scheme to $\PA^-$, the theory of positive parts of ordered semirings, the resulting extension is $\Pi_1$-conservative, as noted in \cite{kaye}, Exercise 7.7.

In this article, we answer Kaye's question in the positive. In fact, following the original conjectural proof strategy, we show that any countable model of $\PA$ can be elementarily extended to an $\omega_1$-like model which carries a full satisfaction class (or, equivalently, which expands to a model of $\CT^-$). We achieve this, in turn, by showing that any countable model of $\CT^-$ satisfying an additional axiom of the internal induction has an end-extension. The proof of this fact, a technical crux of our work, is based on techniques from \cite{lelyk_wcislo_localcollection}, where a partial answer to Kaye's question has been provided, combined with the copying technique from a novel construction of satisfaction classes presented by Pakhomov in his note \cite{pakhomov_copying_satisfaction}.\footnote{The note has not been published, but the ideas contained there were discussed in the presentation \cite{pakhomov_talk_yet_another_satisfaction} available at the Internet address given in the references.}

\section{Preliminaries}

This article concerns truth theories over Peano Arithmetic, $\PA$. Truth theories result by adding to arithmetic a fresh predicate $T(x)$ with the intended reading ``$x$ is the G\"odel code of a true (arithmetical) sentence'' and axioms guaranteeing that $T$ actually displays truth-like behaviour. Crucially for our purposes, Peano arithmetic can formalise syntax, so actually postulating the existence of such a truth predicate makes sense. A comprehensive introduction to truth theories may be found in \cite{halbach} and an extensive treatment of the Tarski Boundary problem may be found in \cite{cies_ksiazka}. 

\subsection{Truth theories}

Let us introduce the main theories of our interest.

\begin{definition} \label{def_ctminus}
	By $\CT^-$ (Compositional Truth) we mean the theory in the arithmetical language $\LPA$ extended with one unary predicate $T(x)$ whose axioms are the axioms of $\PA$ along with the following compositional clauses:
	\begin{enumerate}
		\item $\forall x \Big(T(x) \rightarrow \Sent_{\LPA}(x)\Big).$
		\item $\forall s,t \in \ClTerm_{\LPA}\Big( T(s=t) \equiv \val{s} = \val{t} \Big)$. 
		\item $\forall \phi \in \Sent_{\LPA} \Big(T \neg \phi \equiv \neg T \phi\Big).$
		\item $\forall \phi, \psi \in \Sent_{\LPA} \Big(T (\phi \wedge \psi) \equiv T \phi \wedge T \psi\Big).$
		\item $\forall \phi \in \form_{\LPA}^{\leq 1} \Big(T \exists v \phi(v) \equiv \exists x T \phi(\num{x}) \Big).$
		\item $\forall \bar{s}, \bar{t} \in \ClTermSeq_{\LPA} \forall \phi \in \form_{\LPA} \Big(\overline{\val{s}} = \overline{\val{t}} \rightarrow T \phi(\bar{s}) \equiv T \phi(\bar{t})\Big).$
	\end{enumerate} 
\end{definition}	

The last clause, called \df{the regularity axiom} states that the truth of a sentence does not depend on the exact terms which are used in it, but rather on the values thereof. In the proof, we will actually need a stronger regularity condition which we will discuss in Subsection \ref{ssec_regularity}. 

One of the fundamental facts about $\CT^-$ is that it does not have any arithmetical content extending $\PA$. 

\begin{theorem}[Kotlarski--Krajewski--Lachlan]
	$\CT^-$ is conservative over $\PA$.
\end{theorem}

A number of additional conditions can be imposed on the truth predicate so that the resulting theory remains conservative. One of the most notable is that we can assume that every formula $\phi \in \form_{\LPA}$, considered separately, satisfies induction.

\begin{definition} \label{def_internal induction}
	By the \df{internal induction} axiom ($\INT$) we mean the following statement:
	\begin{displaymath}
		\forall \phi \in \form^{\leq 1}_{\LPA} \Big( T\phi(0) \wedge \forall x \big(T\phi(\num{x}) \rightarrow  T\phi(\num{S(x)})\big) \rightarrow \forall x \phi(\num{x})\Big).
	\end{displaymath}
\end{definition}

 By compositional axioms this is equivalent to saying that all induction axioms are true. Perhaps somewhat surprisingly, this also yields a conservative extension of $\PA$. This result was originally announced in \cite{kkl}. A proof, by a different, model theoretic methods announced in \cite{enayatvisser2}, can be found in the longer, unpublished, and privately circulated manuscript \cite{enayat_visser_long}. Another argument, now purely proof-theoretic, has been presented in
 \cite{leigh} (where it is proved that one can extend $\CT^-$ with an arbitrary statement of the form ``all instances of the axioms scheme $\Gamma$ are true'' while still keeping the theory in question conservative. It is easy to check that the compositional axioms allow us to derive our version of the internal induction from the statement that all the induction axioms are true). 

\begin{theorem}[Kotlarski--Krajewski--Lachlan]
	$\CT^- + \INT$ is conservative over $\PA$. 
\end{theorem}

\subsection{Models of $\PA$}

This article will make use of some classical theory of models of $\PA$. Let us now review some basic facts of this area. The standard references are \cite{kaye} (introductory) and \cite{kossakschmerl} (more advanced), where the proofs of the theorems stated here, and much more, can be found. The first result which we will use repeatedly is resplendence of recursively saturated models. 

\begin{definition} \label{def_resplendence}
	We say that a model $M$ is \df{resplendent} if for any second-order formula $\phi(X)$ with a single second-order variable with all quantifiers ranging over first-order variables and possibly with first-order parametres, there exists $A \subseteq M$ such that $(M,A) \models \phi(A)$. 
\end{definition}

The notion of resplendence is very rich in consequences and yet, resplendent models of strong theories are rather easy to find in nature.

\begin{theorem}[Barwise--Schlipf, Ressayre] \label{tw_barwisschlipf}
	Every countable recursively saturated model of $\PA$ is resplendent.
\end{theorem}

End-extensions of models of truth theories play a crucial role in our paper. They are also a very classical thread in the theory of the models of $\PA$.

\begin{definition}
	Let $M \subseteq N$ be models of $\PA$. We say that $N$ is an \df{end-extension} of $M$ iff for any $c \in N \setminus M$ and an arbitrary $a \in M$, $N \models a<c$. We denote this by $M \subseteq_e N$ (or $M \preceq_e N$ if this is in fact an elementary extension). If $M \neq N$, we call the extension \df{proper}.
\end{definition}

We will also use a more sophisticated variant of extensions. 

\begin{definition} \label{def_coservative_extension}
	Let $M \preceq_e N$. We say that $N$ is a \df{conservative} extension of $M$ iff for any formula $\phi$, possibly with parametres from $N$, there exists a formula $\psi$, possibly with parametres from $M$ such that for any $a \in M$,
	\begin{displaymath}
		N \models \phi(a) \textnormal{ iff } M \models \psi(a).
	\end{displaymath}
\end{definition}

The presence of parametres from $N$ is a crucial requirement in the above definition. Without them, the conclusion would follow trivially by elementarity. One of the crucial facts from the model theory of $\PA$ states that conservative-extensions exist (in fact, they automatically have to be end-extensions).

\begin{theorem}[MacDowell--Specker; Gaifman] \label{tw_macdowellspecker}
	Let $M$ be an arbitrary model over a countable signature satisfying the induction scheme for the full language and all axioms of $\PA$. Then $M$ has a proper elementary conservative end-extension $N$.
\end{theorem}

The statement and the proof of the above result can be found in \cite{kossakschmerl}, Theorem 2.2.8.\footnote{Note that the authors use the following convention: ``model'' means precisely ``a model of a theory over a countable signature satisfying the induction scheme for the full language and all axioms of $\PA$'' as explained on p.1 of the cited source.}

An important class of models are $\kappa$-like models, of which we will need a special case:
\begin{definition} \label{def_omega1}
	Let $M \models \PA$. We say that $M$ is an $\omega_1$-like model if $|M| = \aleph_1$, but any proper initial segment of $M$ is countable.
\end{definition}

The proof of existence of $\omega_1$-like models is the prototype for the argument presented in this work.

\begin{theorem}
	Let $M \models \PA$  be an arbitrary countable model. Then there exists an elementary $\omega_1$-like end-extension $M \prec_e N$.
\end{theorem}
\begin{proof}
	By repeatedly using Theorem \ref{tw_macdowellspecker}, we can construct a sequence of countable models $M_{\alpha}, \alpha<\omega_1$ such that for any $\alpha<\beta$, 
	\begin{displaymath}
		M_{\alpha} \prec_e M_{\beta}.
	\end{displaymath}
	Let $N = \bigcup_{\alpha < \omega_1} M_{\alpha}$. Then $N \models \PA$ as a union of an elementary chain, $N$ has cardinality $\aleph_1$, and for an arbitrary proper initial segment $I$, there exists an element $a \in N$, such that $a \notin I$. Let $\alpha < \omega_1 $ be any ordinal such that $a \in M_{\alpha}$. Then, since $N$ is an end-extension of $M_{\alpha}$,  $I \subset M_{\alpha}$ and hence it is countable. This shows that $N$ is $\omega_1$-like.
\end{proof}

In our paper, we will use a technical condition on cuts of models of $\PA$.

\begin{definition} \label{def_cut}
	Let $M \models \PA$. We say that a subset $I$ of $M$ is a \df{cut} of $M$ if the following conditions hold:
	\begin{itemize}
		\item For any $a \in I$ and any $b \leq a$, $b \in M$.
		\item For any $a \in I$, $S(a) \in I$.
	\end{itemize}
\end{definition}

Cuts in models of arithmetic are a classical subject of intensive study in which a number of their structural properties were isolated. In our argument, we will isolate what we think is a new, very weak, regularity condition, whose definition we postpone until Section \ref{sec_end_ext}. 

\subsection{Truth predicates and satisfaction classes} \label{ssec_truth and satisfaction}

In the literature, there are two competing treatments of the truth-like notions. The first one, more traditional, originates in the theory of models of $\PA$ and speaks of satisfaction classes which are treated primarily as subsets in models of arithmetic. The other stems from discussions in philosophical logic and speaks of truth theories which are treated primarily in the axiomatic manner.

Intuitively, a satisfaction class in a model $M \models \PA$, should be a set $S$ of pairs $(\phi,\alpha)$, where $\phi$ is an arithmetical formula in the sense of the model and $\alpha$ is a $\phi$-assignment such that, collectively, the pairs in $S$ satisfy Tarski's compositional conditions. However, we are often interested in such classes where $S$ is only required to work for some subset of formulae in the model. This, however, makes the notion of a satisfaction class somewhat subtle. In particular, the exact definition of a satisfaction class in not quite consistent between different authors. Below, we present the version from \cite{bw_automorphisms_definability} which, we believe, captures various uses of the notion most smoothly.

\begin{definition} \label{def_satisfaction_class}
	Let $M \models \PA$, let $S \subset M^2$, and let $\phi \in \form_{\LPA}(M)$. We say that $S$ is \df{compositional} at $\phi$ if for any $\alpha \in \Asn(\phi)$, $(\phi,\alpha) \in S$ iff one of the following conditions hold:
	\begin{itemize}
		\item There exist $s,t \in \Term_{\LPA}(M)$ such that $\phi = (s=t)$ and $s^{\alpha} = t^{\alpha}$. 
		\item There exists $\psi \in \form_{\LPA}(M)$ such that $\phi = \neg \psi$ and $(\psi,\alpha) \notin \phi$. 
		\item There exist $\psi, \eta \in \form_{\LPA}(M)$ such that $\phi = \psi \vee \eta$ and $(\psi,\alpha) \in S$ or $(\eta,\alpha) \in S$. 
		\item There exist $\psi, \eta \in \form_{\LPA}(M)$ such that $\phi = \psi \wedge \eta$ and both $(\psi,\alpha) \in S$ and $(\eta,\alpha) \in S$. 
		\item There exists $\psi \in \form_{\LPA}(M)$ and $v \in \Var$ such that $\phi = \exists v \psi$ and $S(\phi,\alpha)$ holds iff $S(\psi,\beta)$ holds for some $\beta \in \Asn(\psi)$ such that $\beta(w)$ is equal to $\alpha(w)$ for all $w$ different from $v$ (which is not required to be in the domain of $\beta$).
		\item There exists $\psi \in \form_{\LPA}(M)$ and $v \in \Var$ such that $\phi = \forall v \psi$ and $S(\phi,\alpha)$ holds iff $S(\psi,\beta)$ holds for all $\beta \in \Asn(\psi)$ such that $\beta(w)$ is equal to $\alpha(w)$ for all $w$ different from $v$.
	\end{itemize}
	If $\alpha$ and $\beta$ are two valuations which are equal except, possibly, on a single variable $v$, we denote this with $\alpha \sim_v \beta$.

	 We say that a set $S \subset M^2$ is a \df{satisfaction class} if there exists a set $D \subset \form_{\LPA}(M)$ such that
	\begin{itemize}
		\item $D$ is closed under taking direct subformulae.
		\item For any $\phi \in D$, $S$ is compositional at $\phi$. 
		\item For any $\phi \in D$ and any $\alpha \in \Asn(\phi)$, either $(\phi,\alpha) \in S$ or $(\neg \phi, \alpha) \in S$. 
		\item For any $\phi \in \form_{\LPA}$, if there exists an $\alpha$ such that $(\phi,\alpha) \in S$ or $(\neg \phi, \alpha) \in S$, then $\phi \in D$. 
	\end{itemize}
	We call the maximal set $D$ satisfying the above conditions a \df{domain} of $S$, denoted $\dom(S)$. If the domain $D$ of $S$ is the whole model $M$, we say that $S$ is a \df{full} satisfaction class. 
\end{definition} 

Unfortunately, it turns out that in the absence of some form of induction, the connection between truth predicates and satisfaction classes is not as clear-cut as one could hope. A discussion of that phenomenon may be found in \cite{bw_automorphisms_definability}. However, the distinction between truth and satisfaction classes trivialises if we assume that certain regularity properties holds and the truth value of a formula does not depend on the choice of specific terms, but rather on their values and it does not depend on the choice of the specific bound variables. We will present the exact assumptions on truth and satisfaction classes in question in the next subsection.

\subsection{Syntactic templates and regularity} \label{ssec_regularity}

As we mentioned before, in the main proof we will need a technical condition on the regularity of the constructed truth and satisfaction classes. We will now spell out those, admittedly tedious, technicalities. Similar considerations played the same role in our previous works, such as \cite{lelyk_wcislo_localcollection} or (in a somewhat different formulation) in \cite{bw_automorphisms_definability}.

\begin{definition} \label{def_syntactic_template}
	Let $\phi$ be an arbitrary formula. By the \df{syntactic template} of $\phi$, we mean the unique formula $\widehat{\phi}$ such that:
	\begin{enumerate}
		\item $\widehat{\phi}$ differs from $\phi$ only by term substitution and renaming bound variables.
		\item There are no complex terms in $\widehat{\phi}$ which contain only free variables. 
		\item The formula $\widehat{\phi}$ contains no closed terms. 
		\item Every free variable occurs in $\widehat{\phi}$ at most once.
		\item Every bound variable in $\widehat{\phi}$ is quantified over only once.
		\item The choice of free and bound variables is such that $\widehat{\phi}$ is the minimal formula satisfying the above conditions.  
	\end{enumerate}
If $\widehat{\phi} = \widehat{\psi}$, we say that $\phi$ and $\psi$ are \df{syntactically similar}. We denote this relation by $\phi \sim \psi$. 
\end{definition}

In essence, a syntactic template represents the pure syntactic tree of a formula in which all the terms which involve no bound variables were erased and replaced by single variables. 

\begin{example}
	Let
	\begin{displaymath}
		\phi = \exists x \forall y  \big(x+(y \times 0) = S(0) + (x \times (z \times v))\big)
	\end{displaymath}
	 Then
	\begin{displaymath}
		\widehat{\phi} = \exists v_0 \forall v_1  \big(v_0+(v_1 \times w_0) = w_1 + (v_0 \times w_2)\big),
	\end{displaymath} 
where the variables $v_i, w_i$ are chosen so that the resulting formula is minimal.
\end{example}

One of the reasons to introduce syntactic templates is because, as already mentioned, we would like to work with satisfaction classes and the connection between satisfaction and truth classes is not quite as neat, as one could expect and some regularity assumptions seem to be required to actually ensure that actually the two notions coincide. Below, $\phi[\alpha]$ is the sentence obtained by substituting the numeral $\num{\alpha(v)}$ for every instance of the variable $v$ in the formula $\phi$. 

\begin{definition} \label{def_syntactic_similarity}
	Let $M \models \PA$. Let $\phi, \psi \in \form_{\LPA}(M)$, let $\alpha \in \Asn(\phi)$, and let $\beta \in \Asn(\psi)$. We say that the pairs $(\phi,\alpha)$ and $(\psi,\beta)$ are \df{syntactically similar} iff $\phi$ is syntactically similar to $\psi$ and there exist sequences of closed terms $\bar{s},\bar{t} \in \ClTermSeq_{\LPA}(M)$ such that:
	\begin{itemize}
		\item  $\bar{\val{s}} =  \bar{\val{t}}$ (the terms in both sequences have the same values);
		\item $\widehat{\phi}(\bar{s})$ differs from $\phi[\alpha]$ only by renaming bound variables.
		\item $\widehat{\phi}(\bar{s}) (= \widehat{\psi}(\bar{s}))$ differs from $\psi[\beta]$ only be renaming bound variables. 
	\end{itemize}
If the pairs $(\phi,\alpha), (\psi,\beta)$ are syntactically similar, we denote this fact with
\begin{displaymath}
	(\phi,\alpha) \sim (\psi,\beta).
\end{displaymath}

\end{definition}

The above notion is actually much simpler than the definition suggests.

\begin{example}
	Let 
	\begin{eqnarray*}
		\phi & : = & \exists x \forall y \ \big( x+(y \times S(S(0))) = z \times (S(0) + S(0)) \big) \\
		\psi & : = & \exists z \forall v \ \big(  z + (v \times (u+w)) = S(S(S(S(0)))) \big).
	\end{eqnarray*}
	Let $\alpha \in \Asn(\phi)$ be an assignment such that
	\begin{displaymath}
		\alpha(z) = 2.
	\end{displaymath}
	Let $\beta \in \Asn(\psi)$ be an assignment such that 
	\begin{displaymath}
		\beta(u) = 1, \beta(w) = 1. 
	\end{displaymath}
Then $(\phi,\alpha) \sim (\psi,\beta)$. Indeed, this is witnessed by $\widehat{\phi} = \widehat{\psi}$ equal to:
\begin{displaymath}
	\exists v_0 \forall v_1 \ \big(v_0 + v_1 \times w_0  = w_1\big). 
\end{displaymath}
and sequences $\bar{s}, \bar{t}$ such that:
\begin{displaymath}
	\bar{s} = \tuple{S(S(0)), S(S(0)) + (S(0) + S(0))}
\end{displaymath}
\begin{displaymath}
	\bar{t} = \tuple{S(0) + S(0), S(S(S(S(0))))}.
\end{displaymath}
Then
\begin{eqnarray*}
	\widehat{\phi}(\bar{s}) & = & \exists v_0 \forall v_1 \ \big(v_0  + (v_1\times S(S(0)))  =  S(S(0)) + (S(0) + S(0))\big) \\
	\phi[\alpha ]& = & \exists x \forall y \ \big(x + (y\times S(S(0))  =  S(S(0)) + (S(0) + S(0))\big)
\end{eqnarray*}
which differ only by substituting closed terms with the same values and renaming variables.
\end{example}

As we already mentioned, we want to restrict our attention to classes for which good regularity properties hold.

\begin{definition} \label{def_syntactic_regularity_}	
	Let $(M,S)$ be a satisfaction class. We say that $S$ is \df{syntactically regular} if for any $\phi,\psi \in \form_{\LPA}(M)$ and $\alpha \in \Asn(\phi), \beta \in \Asn(\psi)$ if $(\phi, \alpha) \sim (\psi,\beta)$, then
	\begin{displaymath}
		(\phi,\alpha) \in S \textnormal{ iff } (\psi,\beta) \in S.
	\end{displaymath}

	Let $(M,T) \models \CT^-$. We say that $T$ is \df{syntactically regular} if for any $\phi, \psi  \in \Sent_{\LPA}(M)$ if $(\phi,\emptyset) \sim (\psi, \emptyset)$ (note that $\emptyset$ is a trivial valuation), 
	\begin{displaymath}
		\phi \in T \textnormal{ iff } \psi \in T.
	\end{displaymath}	
\end{definition}

Under these regularity assumptions, we can actually make a straightforward connection between truth predicates and satisfaction classes.

\begin{proposition} \label{stw_satisfaction_iff_truth}
	Let $M \models \PA$ and let $S$ be a full syntactically regular satisfaction class on $M$. Let 
	\begin{displaymath}
		T = \set{\phi \in \Sent_{\LPA}(M)}{(\phi,\emptyset) \in S}.
	\end{displaymath}
	Then $(M,T) \models \CT^-$ and, moreover, $T$ is syntactically regular. 
	
	Conversely, suppose that $(M,T) \models \CT^-$ and that $T$ is syntactically regular. Let 
	\begin{displaymath}
		S = \set{(\phi,\alpha) \in M^2}{\phi \in \form_{\LPA}(M), \alpha \in \Asn(\phi) \textnormal{, and } \phi[\alpha] \in T}.
	\end{displaymath}
Then $S$ is a full syntactically regular satisfaction class.
\end{proposition}

 If $S$ and $T$ are interdefinable in the way postulated in the above proposition, we say that $S$ is a satisfaction class \df{corresponding} to $T$ and that $T$ is a truth predicate \df{corresponding} to $S$.
A slightly different statement, in which we used weaker regularity assumptions appeared as Proposition 15 in \cite{bw_automorphisms_definability}. The definition of syntactic regularity formulated in this work assumes that a syntactically regular satisfaction class is closed under renaming bound variables (equivalently, under $\alpha$-conversion). This is not needed to obtain the above correspondence, but will be used in the proof of the main theorem. However, \emph{some} regularity assumptions apparently \emph{are} needed in order to obtain the above simple correspondence as discussed in the cited article.

Crucially for our paper, the syntactic regularity conditions can be added as a an additional requirement essentially to any reasonable conservative theory of truth. In particular, we have:

\begin{proposition} \label{stw_consiervativity_of_int_plus_regularity}
	Let $M \models \PA$ be a countable recursively saturated model. Then there exists $T \subset M$ such that $(M,T) \models \CT^- + \INT$ and $T$ is syntactically regular.
\end{proposition}
This fact appeared as Theorem 23 in \cite{lelyk_wcislo_localcollection}. The proof of this fact was strictly speaking omitted, but the proof is a (completely straightforward) modification of an argument which appeared there with a precise comment on what (also straightforward) modification is needed.

\section{The conservativity of collection} \label{sec_cons_coll}

In this section, we present the strategy for the proof of conservativity of the compositional truth predicate with the collection axioms, relegating the demonstration of the crucial technical results to next sections. Let us start with a basic observation already suggested by Kaye as the main tool for the argument we present in this article. 

\begin{proposition} \label{stw_omega_1_ma_kolekcje}
	Suppose that $M \models \PA$ is an $\omega_1$-like model and that $S \subset M$ is an arbitrary subset. Then the expansion $(M,S)$ satisfies the full collection scheme.
\end{proposition}

\begin{proof}
	Let $\phi$ be an arbitrary formula in the language $\LPA$ expanded with the symbol $S$. Fix any $a \in M$ and suppose that 
	\begin{displaymath}
		(M,S) \models \forall x < a \exists y \ \phi(x,y).
	\end{displaymath}	
	Let $f: M \to M$ be a function such that for any $x \in M$ which is smaller than $a$, 
	\begin{displaymath}
		(M,S) \models \phi(x,f(x)).
	\end{displaymath}
	Since $M$ is $\omega_1$-like model, the interval $[0,a]$ is countable. Therefore, the image $f\big[[0,a]\big]$ cannot be cofinal in $M$. This means that there exists $b \in M$ such that  $f\big[[0,a]\big] \subseteq [0,b]$. In particular,
	\begin{displaymath}
		(M,S) \models \forall x <a \exists y< b \ \phi(x,y).
	\end{displaymath} 
	Since $\phi$ was arbitrary, $(M,S)$ is a model of the full collection.
\end{proof}

As a matter of fact, our argument for conservativity of collection follows exactly the path suggested by the above result. 

\begin{theorem} \label{tw_omega_1_like_models}
	Let $M \models \PA$ be an arbitrary countable model. Then, there exists an $\omega_1$-like elementary extension $M' \succ M$ and a full syntactically regular satisfaction class $S \subset M'$. 
\end{theorem}

The main goal of the article is to prove Theorem \ref{tw_omega_1_like_models}. One obvious potential strategy for a proof would be to show that for an arbitrary countable model $(M,T) \models \CT^-$, one can find a proper end-extension $(M',T') \succ_e (M,T)$. However, this direct strategy cannot quite work, since  as demonstrated by \cite{smith}, countable models of $\CT^-$ do not necessarily have end-extensions.

\begin{theorem}[Smith] \label{tw_model_ctminus_bezrozszerzen}
	For any countable recursively saturated model $M \models \PA$, there exists $T \subset M$ such that $(M,T) \models \CT^-$ and there is no $(M',T') \succ_e (M,T)$ satisfying $\CT^-$.  
\end{theorem}

The proof uses the following result which can be found in \cite{smith},Theorem 3.3.

\begin{theorem}[Smith] \label{tw_formuly_w_ctminus_chwytaja_zbiory_rek_nas}
	Let $M$ be a countable recursively saturated model of $\PA$. Let $A \subseteq M$ be an arbitrary set such that the expansion $(M,A)$ is recursively saturated in the extended language. Then there exist $T \subset M$ and $\phi(v) \in \form^{\leq 1}_{\LPA}$ such that $(M,T) \models \CT^-$ and 
	\begin{displaymath}
		A = \set{x \in M}{(M,T) \models T\phi(\num{x})}.
	\end{displaymath}
\end{theorem}
By resplendence of the countable recursively saturated models, we obtain the following corollary.
\begin{corollary} \label{cor_formuly_w_ctminus_chwytaja_zbiory_o_dowolnych_teoriach}
	Let $M$ be a countable recursively saturated model of $\PA$, let $a_1, \ldots, a_n \in M$, and let $A \subseteq M$ be an arbitrary set. $T \subset M$ and $\phi(v) \in \form^{\leq 1}_{\LPA}$ such that $(M,T) \models \CT^-$ and
	\begin{displaymath}
		(M,A) \equiv (M,A'),
	\end{displaymath} 
	where
	\begin{displaymath}
		A' = \set{x \in M}{(M,T) \models T\phi(\num{x})}.
	\end{displaymath} 
\end{corollary}

Now, we can prove theorem \ref{tw_model_ctminus_bezrozszerzen}.
\begin{proof}
	Let $M \models \PA$ be countable and recursively saturated. Pick any $a \in M$ in the nonstandard part. Then there exists a bijection between $[0,a]$ and $M$. By Corollary \ref{cor_formuly_w_ctminus_chwytaja_zbiory_o_dowolnych_teoriach}, there exist an expansion $(M,T) \models \CT^-$ and a formula $\phi(x,y) \in \form_{\LPA}(M)$ such that the following hold:
	\begin{itemize}
		\item $(M,T) \models \forall x < a \exists! y \ T \phi(\num{x},\num{y})$.
		\item $(M,T) \models \forall y \exists! x < a \ T \phi(\num{x},\num{y})$.  
	\end{itemize} 
Notice that by the compositional axioms, these conditions can be equivalently rewritten as:
	\begin{itemize}
		\item $(M,T) \models T \big( \forall x < \num{a} \exists! y \phi(x,y) \big)$.
		\item $(M,T) \models T \big( \forall y \exists! x < \num{a} \ T \phi(x,y) \big)$.  
	\end{itemize} 

Now, suppose that $(M',T') \supset_e (M,T)$ satisfies $\CT^-$. Then in $(M',T')$, $\phi$ defines a bijection between $[0,a]^{M'}$ and $M'$ under the truth predicate. Moreover, this bijection  extends the one defined in $(M,T)$. However this is impossible, since $[0,a]^{M'} = [0,a]^{M}$, but $M \neq M'$.
\end{proof}

Notice that in the above example, we used the fact that collection was violated for a formula of the form $T\phi(\num{x},\num{y})$. One can wonder whether this is in fact the only possible obstruction to the existence of end-extensions. It turns out that the answer to this question is indeed positive. For simplicity, we will first prove the end-extension result under slightly stronger assumptions, which we will relax in Section \ref{sec_internal_collection}.

\begin{theorem} \label{tw_modele_ctminus_plus_int_maja_rozsszerzenia}
	Let $(M,T)  \models \CT^- + \INT$ be a countable model with $T$ syntactically regular. Then there exists a proper end-extension $(M,T) \subset_e (M',T') \models \CT^- + \INT$ satisfying the syntactic regularity condition.
\end{theorem}

The proof of the above result is the heart of this work. We will present it in the next section. Before we do, let us draw our main corollary. 

\begin{theorem} \label{tw_omega1_modele_ctminus_istnieja}
	Let $M \models \PA$ be an arbitrary countable model. Then there exists an elementary extension $M \prec M'$ such that $M'$ is $\omega_1$-like and there exists $T \subset M'$ for which $(M',T) \models \CT^-$. 
\end{theorem}

\begin{proof}
	Fix an arbitrary countable model $M \models \PA$. Let $M_0 \succeq M$ be a countable, recursively saturated model of $\PA$. By resplendence, we can find an expansion of $M_0$ to a model $(M_0,T_0) \models \CT^- + \INT$ with $T_0$ syntactically regular. We define inductively a sequence of countable models $(M_{\alpha},T_{\alpha}), \alpha < \omega_1$ satisfying $\CT^- + \INT$ as follows:
	
	\begin{itemize}
		\item $(M_{\alpha+1},T_{\alpha+1}) \supset_e (M_{\alpha},T_{\alpha})$ is an arbitrary proper end-extension satisfying $\CT^- + \INT$. 
		\item $(M_{\lambda},T_{\lambda}) = \bigcup_{\gamma < \lambda} (M_{\gamma}, T_{\gamma})$ for limit $\lambda$. 
	\end{itemize}

	The models at successor steps can constructed by Theorem \ref{tw_modele_ctminus_plus_int_maja_rozsszerzenia}, so we only have to check that the induction hypotheses can be maintained in the limit steps. However, it can be checked in a straightforward manner that the compositional clauses of $\CT^-$ are preserved in the unions of models. Also notice that over $\CT^-$, the internal induction is equivalent to a $\Pi_1$-sentence saying ``all the instances of the arithmetical induction scheme are true'', so it is also preserved at the limit steps. 
\end{proof}

We can now complete the main line of the argument.

\begin{theorem} \label{tw_coll_is_conservative}
	$\CT^- + \Coll$ is conservative over $\PA$. 
\end{theorem}
\begin{proof}
	It is enough to show that for an arbitrary arithmetical sentence $\phi$, if $\PA + \phi$ is consistent, then $\CT^- + \Coll + \phi$ is consistent.
	
	Suppose that $\PA + \phi$ is consistent and take any countable model $M \models \PA + \phi$. By Theorem \ref{tw_omega1_modele_ctminus_istnieja}, there exists an elementary $\omega_1$-like extension $M' \succeq M$ with a subset $T \subset M'$ such that $(M',T) \models \CT^-$. By Proposition \ref{stw_omega_1_ma_kolekcje}, $(M',T)$ is actually a model of $\CT^- + \Coll$. By elementarity, $M' \models \phi$, so $(M',T)$ witnesses the consistency of $\CT^- + \Coll + \phi$.  
\end{proof}

Actually, since our methods really rely on the countability of the models involved in the construction, it is unclear whether we can show the existence of $\kappa$-like models of $\CT^-$ for an arbitrary $\kappa$.

\begin{question}
	Let $M \models \PA$ be an arbitrary countable model and let $\kappa$ be an arbitrary cardinal. Does there exist a $\kappa$-like model $M' \succ_e M$ which expands to a model of $\CT^-$?
\end{question}

\section{End-extensions of satisfaction classes} \label{sec_end_ext} 

In this section, we will prove the main result. As we mentioned in the previous section, the key step is to prove that the countable models of $\CT^- + \INT$ always have end-extensions. For the technical convenience, in this section we will switch to the language of regular satisfaction classes. By Proposition \ref{stw_satisfaction_iff_truth}, in our framework they directly correspond to the models of $\CT^-$ satisfying full regularity.

The structure of argument will be divided into two main parts:
\begin{itemize}
	\item First we show that given a countable model $M$ of $\PA$ carrying a full regular satisfaction class which satisfies the internal induction, we can end-extend it to a model $M'$ with a partial satisfaction class whose domain includes all formulae with the syntactic depth in $M$ and such that $M$ is a nicely behaved cut of $M'$.  
	\item Then we prove that if a countable model $M'$ has a partial satisfaction class whose domain consists of formulae with the syntactic depth in a certain nicely behaved cut, then in $M'$ we can find a full satisfaction class. Moreover, internal induction can be preserved in this extension.
\end{itemize}

The next two subsections will be devoted to those two main steps of the proof. In particular, wee will make precise  the requirements we impose on the cuts in question. 

\subsection{Stretching Lemma}

In this part, we will discuss the first of the main steps in the end-extension theorem. The key ideas of this subsection appeared already in \cite{lelyk_wcislo_localcollection}. However, since we will need to extract some additional information from the construction, we will present here the full proof.

As we already mentioned, we will need to impose certain regularity conditions on the cuts arising in our construction. The exact choice of the condition is rather subtle.

\begin{definition} \label{def_local_semiregularity}
	Let $I \subset M \models \PA$ be a nonstandard cut. We say that $I$ is \df{locally semiregular} in $M$ if for any $a \in I$ and any function $f: [0,a] \to M$, there exist a nonstandard $a' \leq a $ and $b \in I$ such that the following condition holds:
	\begin{displaymath}
		f[[0,a']] \cap I \subseteq [0,b].
	\end{displaymath}
\end{definition} 

A word of comment is certainly in place. One of the conditions on cuts, classically investigated in the theory of models of $\PA$, is semiregularity. We say that a cut $I$ is \df{semiregular} in $M$ if for any function $f$ coded in $M$ with the domain $[0,a]$ for some $a \in I$, the set of thhe values $f(i)$ such that $f(i) \in I$ is bounded in $I$. Local semiregularity wekens this condition, demanding instead that the function $f$ can be restricted to a nonstandard initial segment so that the bounding condition holds. 

Admittedly, this is a somewhat technical requirement. However, the choice of this exact condition will be crucial in Section \ref{sec_internal_collection}, since we were not able to show the main end-extension result from that part which would guarantee  any stronger regularity requirements on cuts. On the other hand, any weaker conditions known in the literature do not seem to suffice to perform the copying construction of Theorem \ref{tw_copying}.

Let also us add that we were also not able to find any condition previously known in the literature which would be equivalent to local semiregularity.

\begin{lemma}[Stretching Lemma] \label{lem_stretching}
	Let $M \models \PA$ be a countable model and let $S \subset M^2$ be a full regular satisfaction class satisfying the internal induction. Then there exists a proper end extension $(M,S) \subset (M',S')$ such that:
	\begin{itemize}
		\item $S'$ is a regular partial satisfaction class.
		\item For any $\phi \in \form_{\LPA}(M')$, $\phi \in \dom(S')$ iff $\dpt(\phi) \in M$. 
		\item $S'$ satisfies the internal induction axiom. 
		\item The extension $M \prec_e M'$ is elementary and $M$ is a locally semiregular cut of $M'$. 
	\end{itemize}
\end{lemma}

\begin{proof}
	Fix a model $(M,S)$ as above. For any $\phi \in \form_{\LPA}(M)$, let $S_{\phi}$ be the set:
	\begin{displaymath}
		S_{\phi} = \set{\alpha}{(\phi,\alpha) \in S}.
	\end{displaymath}

	Consider the model $(M,S_{\phi})_{\phi \in M}$. By the internal induction axiom, this model satisfies full induction. Since $M$ is countable, the signature of the expanded model has only countably many symbols. Thus by MacDowell-Specker Theorem, there exists a proper conservative elementary end-extension
	\begin{displaymath}
		(M,S_{\phi})_{\phi \in M} \prec_e (M',S'_{\phi})_{\phi \in M}.
	\end{displaymath}
	
	Now, let $S' \subset (M')^2$ be defined with the condition: $(\phi,\alpha) \in S'$ iff there exists a pair $(\psi,\beta)$ such that:
	\begin{itemize}
		\item $(\psi,\beta) \sim (\phi,\alpha)$.
		\item $\psi \in M$.
		\item $S'_{\psi}(\beta)$ holds.
	\end{itemize}

In other words, we take a union of sets $S'_{\psi}$ and close it under syntactic similarity. We can then check that by elementarity of the extension $(M',S'_{\phi})_{\phi \in M}$ and by regularity of $S$, the constructed set $S'$ is a regular satisfaction class satisfying internal induction. The regularity of $S'$ follows directly by construction.

We have to check that $M$ is a locally semiregular cut of $M'$. In fact, we will check that it is semiregular. Let 
\begin{displaymath}
	f: [0,a] \to M'
\end{displaymath}
be a function coded in $M'$. We want to check that for any $b \in M$, there exists $c \in M$ such that the values of $f \res M$ are bounded by $c$. Since $f$ is coded in  $M'$, by conservativity of the extension $(M,S_{\phi})_{\phi \in M} \prec (M',S'_{\phi})_{\phi \in M}$, the set $f \cap M$ is definable in the former structure. However, since this model satisfies full induction in the expanded language, the values of $(f \res [0,b]) \cap M$ have to be bounded in $M$. 
\end{proof}

We will discuss in Section \ref{sec_internal_collection} how internal induction can be eliminated from the above argument.

\subsection{Copying Lemma}

In the previous subsection, we have shown how to extend a satisfaction class ``upwards'' so that it is still defined for all the formulae in the original model. In this subsection, we will discuss the essentially novel part of our argument: we will show how, under additional model-theoretic assumptions, we can extend a satisfaction class defined on a cut of formulae to the whole model. Our argument crucially uses ideas introduced by Fedor Pakhomov in his construction of a satisfaction class presented in his unpublished note \cite{pakhomov_copying_satisfaction}.\footnote{The construction described in the note was discussed in the Autumn 2020 in the talk \cite{pakhomov_talk_yet_another_satisfaction} seminar ``Epistemic and Semantic Commitments of Foundational Theories.'' The slides and the recording of the talk are available at the seminar webpage given in the bibliography.} 

\begin{theorem}[Copying Lemma] \label{tw_copying}
	Let $M \models \PA$ be a countable model. Let $I \subset M$ be a locally semiregular cut in $M$. Suppose that there exists a syntactically regular satisfaction class $S \subset M^2$ whose domain consists of formulae with depth in $I$. Then there exists a full regular satisfaction class $S' \supset S$. Moreover, if $S$ satisfies internal induction, then $S'$ also satisfies it. 
\end{theorem}

\begin{proof}
	Let, $M,I,S$ be as in the assumptions of the theorem. Let $D$ be the set of all formulae in $M$ whose syntactic depth is in $I$ (i.e., $D$ is the domain of $S$). 
	We will construct a function $f: \form_{\LPA}(M) \to D$ satisfying the following conditions for any $\phi, \psi \in \form_{\LPA}$ and any $v \in \Var(M)$:
	\begin{itemize}
		\item $f \res D = \id$. 
		\item For any $\phi, \psi$ if $\phi \sim \psi$, then $f(\phi) \sim f(\psi)$.
		\item $f(\neg \phi) = \neg f(\phi).$
		\item $f(\phi \wedge \psi) = f(\phi) \wedge f(\psi).$
		\item $f(\phi \vee \psi) = f(\phi) \vee f(\psi).$
		\item $f(\exists v \phi) = \exists v f(\phi)$.
		\item $f(\forall v \phi) = \forall v f(\phi)$. 
	\end{itemize}

For a formula $\phi \in \form_{\LPA}(M)$ and $a \in M$, let $U(\phi,a)$ be the set of formulae $\psi$ such that $\psi \sim \psi'$ for some $\psi'$ located at most at height $a$ in the syntactic tree of $\phi$. 

We will define a sequence of nonstandard elements of $M$, $a_0> a_1>a_2 \ldots$ and functions $f_0, f_1, \ldots$ such that each function $f_i$ satisfies the following conditions:
\begin{itemize}
	\item $f_i: U(\phi_i,a_i) \to \form_{\LPA}(M)$. 
	\item If $i \leq k$, then $f_i \res U(\phi_i,a_k) \cap U(\phi_k,a_k) = f_k \res U(\phi_i,a_k) \cap U(\phi_k,a_k)$.
	\item $f_i$ commutes with quantifiers and connectives in the sense postulated for $f$. 
	\item $f_i \res D = \id$.
	\item If $\phi \sim \psi$, then $f_i(\phi) \sim f_i(\psi)$.  
\end{itemize}
Notice that we \emph{do not} require that $f_0 \subset f_1 \subset f_2$ or that the functions $f_i$ agree on the whole common domain. Finally, we will set: 
\begin{displaymath}
	f(\phi) = f_i(\phi),
\end{displaymath}
where $\phi_i = \widehat{\phi}$ (so $i$ is the index of the template of $\phi$ in our enumeration).

Let us check that $f$ defined in this way indeed satisfies our conditions. It is enough to check that $f$ preserves syntactic operations, since the other conditions follow directly by assumption on the functions $f_i$. We will check the condition for conjunction, the others being similar or completely analogous. So let us fix formulae $\phi, \psi$. Suppose that $\phi \sim \phi_k, \psi \sim \phi_l, \phi \wedge \psi \sim \phi_m$. Let $n = \max(k,l,m).$ Let $U_1 = U(\phi_m,a_n) \cap U(\phi_k,a_n), U_2 = U(\phi_m,a_n) \cap U(\phi_l, a_n)$. Then 
\begin{displaymath}
	f_k \res U_1 = f_m \res U_1
\end{displaymath}
and 
\begin{displaymath}
	f_l \res U_2 = f_m \res U_2.
\end{displaymath}
In particular $f(\phi) = f_m(\phi)$, $f(\psi) = f_m(\psi)$ and the claim follows, since $f_m$ preserves the syntactic structure. So it is enough to construct the sequences $(f_i), (a_i)$ as above. We will also construct an auxiliary sequence $(c_i)$ of the elements of $M$.

Let $a_0$ be an arbitrary element of the cut $I$.  The construction of $f_0$ is very similar to the construction of functions $f_{i+!}$ in the successor steps, so we will go directly to that case indicating the (small and obvious) differences whenever they appear in the proof. 

Let $a:=a_{n+1}$ be an arbitrary nonstandard element of $I$ such that:
\begin{displaymath}
	a_{n+1}2^{a_{n+1}} < \frac{a_n}{2}.
\end{displaymath}
Let $\phi:= \phi_{n+1}$. Now, consider the set $U(\phi,a)$. Since $a \in I$, by local semiregularity of this cut we can conclude that (possibly after replacing $a$ with a smaller element $a'$ which we will for simplicity still denote $a$) there exists $c :=c_{n+1} \in I$ such that all subformulae of $\phi_{n+1}$ which occur at the $a$-th level of the syntactic tree of $\phi$ have themselves syntactic depth either $< c$ or not in $I$. 

Consider the following relation $\unlhd$: 
\begin{displaymath}
	\xi \unlhd \eta
\end{displaymath} 
iff there exists a coded sequence of formulae:
\begin{displaymath}
	\xi = \xi_0, \xi_1, \ldots, x_p = \eta
\end{displaymath}
such that for any $i$, $\xi_i$ is a direct subformula of $\xi_{i+1}$ and $\xi_i \in U(\phi,a)$. In other words, $\xi \unlhd \eta$ means that $\xi$ is a subformula of $\eta$ as can be witnessed using only formulae from $U(\phi,a)$. 

We want to define the function $f:=f_{n+1}: U(\phi,a) \to \form_{\LPA}(M)$ so that it commutes with the syntactic operations. We actually will need a small technical definition. Let us say that a formula $\phi \in U(\phi,a)$ is \df{weakly minimal} with respect to $\unlhd$ if it has a direct subformula which is not in $U(\phi,a)$. Such formulae can be different from $\unlhd$-minimal formulae in the case their main operator is binary, one of the formulae is in $U(\phi,a)$, and the other is not. 

  It is enough to define the function $f$ on the set weakly minimal formulae from $U(\phi,a)$ and then extend the definition to the whole $U(\phi,a)$ by induction (applied internally in the model). In this manner, we will obtain a function $f_{n+1}$ commuting with all connectives and quantifiers, whenever a formula and all its direct subformulae are in the domain.

Now pick any weakly $\unlhd$-minimal formula $\psi$. We consider two cases:

\paragraph*{Case I} There exists $k < n+1$ such that $\psi \in \dom(f_k)$. Then we set $f(\psi) = f_k(\psi)$, where $k$ is the maximal such index. (When we construct the function $f_0$, we simply omit this step.)

\paragraph*{Case II} Otherwise, let $f(\psi)$ be the formula $\psi$ with every subformula at the level $c$ in the syntactic tree replaced with the sentence $0=0$. 

As mentioned before, the function $f$ can be then uniquely extended to the set $U(\phi,a)$ by induction on $\unlhd$ performed in the model $M$. So it is enough check that the sequence of functions $(f_n)$ defined above satisfies our requirements. 

It is clear by definition that for all $n$, $f_n$ is defined on $U(\phi_n,a_n)$ and that it preserves the syntactic operations. We need to check that these functions identity when restricted to $D$, that they are congruent with respect to $\sim$, and that they satisfy the agreement condition.

\paragraph*{The identity condition} We prove by induction on $n$ that $f_n \res D = \id$. The initial case will be very similar to the induction step, so we only present the latter. 

Fix any formula $\psi \in D \cap U(\phi_{n+1},a_{n+1})$ and first assume that $\psi$ is $\unlhd$-minimal. If $\psi \in \dom(f_k)$ for some $k \leq n$, then $f_{n+1}(\psi) = f_l(\psi)$ where $l$ is the maximal index for which $\psi \in \dom(f_l)$. Then by induction hypothesis $f_l(\psi) = \psi$ which proves the claim. 

If $\psi \notin \dom(f_k)$ for $k \leq n$, then $f_{n+1}(\psi)$ is the formula $\psi$ with any subformula at the syntactic level $c_{n+1}$ replaced with a sentence $0=0$. However, by construction any formula with the syntactic depth from $I$ which belongs to $U(\phi_{n+1},a_{n+1})$ has syntactic depth strictly less than $c_{n+1}$, so in fact the described substitution is trivial and $f_{n+1}(\psi) = \psi$. Then it is enough to observe that on non-minimal formulae $f_{n+1}$ is defined by induction on $\unlhd$ which clearly preserves the identity condition, since the set $D$ is closed under subformulae. 

\paragraph*{The congruence condition} We check by induction on $n$ that for any $\psi,\eta \in \dom(f_n)$ if $\psi \sim \eta$, then $f_n(\psi) \sim f_n(\eta)$. 

First, suppose that $\psi, \eta$ are weakly minimal in $U(\phi_n,a_n)$. First suppose that $\psi \in \dom(f_k)$ for some $k < n$. Notice that the set $U(\phi_k,a_k)$ is closed under $\sim$, so in such case both $\phi$ and $\psi$ are in $\dom(f_k)$ and by induction hypothesis, $f_n(\psi) \sim f_n(\eta)$. 

If, on the other hand, $\psi, \eta \notin \dom(f_k)$ for any $k < n$, then $f_n(\psi)$ and $f_n(\eta)$ are defined by substituting $0=0$ for any subformula at the syntactic depth $c$ in these formulae. Notice that if $\psi \sim \eta$, then their syntactic trees were equal up to term substitutions and renaming bound variables. In such a case, those trees with $0=0$ substituted for all formulae at the depth $c$ will still be equal. 

Now, we can check by induction that $f_n$ is a congruence with respect to $\sim$ on the whole $U(\phi_n,a_n)$ by induction on $\unlhd$, applied internally. 

\paragraph*{The coherence condition} 

We want to check that if $i \leq k$, then
\begin{displaymath}
	f_i \res U(\phi_i,a_k) \cap U(\phi_k,a_k) = f_k \res U(\phi_i,a_k) \cap U(\phi_k,a_k)
\end{displaymath}
We can inductively assume that the desired equality holds for any $j \leq i$:
\begin{displaymath}
	f_i \res U(\phi_i,a_j) \cap U(\phi_k,a_j) = f_j \res U(\phi_i,a_j) \cap U(\phi_j,a_j).
\end{displaymath} 
 Fix any $\phi \in U(\phi_i,a_k) \cap U(\phi_k,a_k).$ Consider any weakly minimal $\psi \in U(\phi_k,a_k)$ which is $\unlhd$-below $\phi$ in $U(\phi_k,a_k)$. 
 
 Now, observe that any $\unlhd$-chain in $U(\phi_k,a_k)$ has at most 
 \begin{displaymath}
 	a_k2^{a_k}
 \end{displaymath}
elements, since for any $\sim$-representative $\eta$ of $\phi_k$, there are at most that many subformulae of $\eta$ in $U(\phi_k,a_k)$ (since they form a tree with at most binary branches of height at most $a_k$). Since $a_k2^{a_k} + a_k < a_i$, all such weakly minimal formulae $\psi$ are in $U(\phi_i,a_i)$ and, in fact in $U(\phi_i,a_j)$ for any $j<k$. Now, by construction, for any such formula $\psi$ we have:
\begin{displaymath}
	f_k(\psi) = f_m(\psi),
\end{displaymath}
where $m \geq i$ is the greatest index smaller than $k$ such that $\psi  \in U(\phi_m,a_m)$. In particular, by the above remark 
\begin{displaymath}
	\psi \in U(\phi_m,a_m) \cap U(\phi_i,a_m).
\end{displaymath}
Hence, by induction hypothesis, 
\begin{displaymath}
	f_m(\psi) = f_i(\psi),
\end{displaymath}
and, consequently,
\begin{displaymath}
	f_k(\psi) = f_m(\psi) = f_i(\psi). 
\end{displaymath}
Since $f_k(\phi)$ is uniquely determined by the values on weakly minimal formulae which are $\unlhd$-smaller than $\phi$ in $U(\phi_k,a_k)$ and as we have just argued, all such chains are also contained in $U(\phi_i,a_i)$, we conclude that
\begin{displaymath}
	f_k(\phi) = f_i(\phi).
\end{displaymath}
This concludes the proof of the coherence clause and therefore of Theorem \ref{tw_copying}.
\end{proof}

\subsection{The proof of the main theorem}

In this section, we will put together the findings of the two previous parts in order to prove Theorem \ref{tw_modele_ctminus_plus_int_maja_rozsszerzenia}, showing that any countable model of $\CT^- + \INT$ has an end-extension. 

\begin{proof}[Proof of Theorem \ref{tw_modele_ctminus_plus_int_maja_rozsszerzenia}]
	As in the formulation of the theorem, let $(M,T) \models \CT^- + \INT$ be a countable model. Let $S$ be a full regular satisfaction class corresponding to $T$. By Lemma \ref{lem_stretching}, there exists a proper end-extension 
	\begin{displaymath}
		(M,S) \subset_e (M',S')
	\end{displaymath}
	such that $S'$ is a regular partial satisfaction class whose domain consists of formulae with depth in $M$ with $M$ being a locally semiregular cut in $M$. Morevoer, the obtained satisfaction class satisfies the internal induction axiom.
	
	Now by Theorem \ref{tw_copying}, there exists $S'' \supset S'$ such that $(M',S'')$ is a full syntactically regular satisfaction class with the internal induction. Let $T'$ be a truth class corresponding to $S''$. Then $(M',T') \models \CT^- + \INT$ is a desired end-extension.  
\end{proof}

\section{Internal collection} \label{sec_internal_collection}

In the previous sections, we have shown that for any complete consistent theory in $\LPA$ extending Peano Arithmetic, there exist a models $M$ satisfying that theory which features a full truth class satisfying the full collection scheme.

However, the fact that we needed to use the internal induction as a vehicle for obtaining end-extensions seems to be an artefact of our argument rather than genuine necessity, especially since internal induction does not seem to follow from full collection. Therefore, it is natural to ask, whether an analogue of Theorem \ref{tw_modele_ctminus_plus_int_maja_rozsszerzenia} holds for some weaker theories. A natural candidate for such an analogue seems to be the principle of the internal collection:

\begin{definition} \label{def_internal_collection}
	By the \df{internal collection} axiom ($\INTColl$), we mean the following principle:
	\begin{displaymath}
		\forall \phi \in \form_{\LPA}^{\leq 1} \forall a \Big(\forall x<a \exists y \ T\phi(\num{x}, \num{y}) \rightarrow \exists b \forall x<a \exists y<b \ T \phi(\num{x}, \num{y}) \Big).
	\end{displaymath}
\end{definition}

Mimicking the usual proof that induction implies collection, we can show using compositional conditions that the internal induction axiom entails the internal collection axiom. In particular, both are conservative over $\PA$. Most likely, the reverse implication does not obtain.

As we already mentioned, internal induction can indeed be replaced with internal collection in our argument.

\begin{theorem} \label{th_int_coll_has_end_extensions}
	Let $(M,T) \models \CT^- + \INTColl$ be a countable model. Then there exists an end-extension $(M,T) \subset_e (M',T') \models \CT^- + \INTColl$. 
\end{theorem}

The proof of the theorem is entirely parallel to the argument for internal induction. The only place in our argument in which internal induction was actually used was to assure in the proof of Stretching Lemma \ref{lem_stretching} that in the end extension $(M,T) \subset_e (M',T')$ we can arrange $M$ to a locally semiregular cut in $M'$. We actually invoked a significantly stronger fact. By using MacDowell--Specker Theorem, we assured that we produced a \emph{conservative} end-extension of models
\begin{displaymath}
	(M,S_{\phi})_{\phi \in M} \prec_e (M',S'_{\phi})_{\phi \in M'}
\end{displaymath}
 and concluded that $M$ is a semiregular cut in $M'$. Unfortunately, if $(M,S)$ satisfies only internal collection, in general the model $(M,S_{\phi})_{\phi \in M}$ will not satisfy full induction, so MacDowell--Specker theorem cannot be applied. However, it turns out that for models with full collection we can find a weaker statement which nevertheless fully suffices for our proof. 

In the proof of Theorem \ref{th_int_coll_has_end_extensions}, we will use more than just slicing having an extension. The following result by Smith (\cite{smith}, Theorem 3.1) will play a crucial role:

\begin{theorem} \label{th_rec_sat_internally_definable}
	Let $M \models \PA$ and let $S$ be a full satisfaction class on $M$. Let $\phi \in \form_{\LPA}(M)$. Let 
	\begin{displaymath}
		S_{\phi} = \set{\alpha \in \Asn(\phi)}{(\phi,\alpha) \in S}.
	\end{displaymath}
Then the expansion $(M,S_{\phi})$ is recursively saturated.
\end{theorem}

As an immediate corollary, we see that the slicing $(M,S_{\phi})_{\phi \in M}$ of a model with a full satisfaction class satisfies a good deal of saturation. Let us give this kind of saturation a handy name: 
\begin{definition}
	Let $M$ be an arbitrary model. We say that $M$ is \df{piecewise recursively saturated} if any recursive type $p$ in which only finitely many symbols from $M$ are used is realised in $M$.  
\end{definition}
\begin{corollary} \label{cor_slicing_piecewise_saturated}
	Let $M \models$ be a countable model with a full satisfaction class $S$ satisfying the internal collection. Then the slicing $(M,S_{\phi})_{\phi \in M}$ is a piecewise recursively saturated structure in a countable language satisfying the full collection scheme. 
\end{corollary}
\begin{proof}
	This is an immediate corollary to Theorem \ref{th_rec_sat_internally_definable}. If we expand $M$ with any finite set of predicates $S_{\phi_1}, \ldots, S_{\phi_n}$, the resulting structure is recursively saturated since these predicates can be defined from one predicate $S_{\psi}$ using the arithmetical coding of tuples.
\end{proof}
It seems that the use of coding was in fact not necessary in the above argument and we could simply reprove Smith's result directly for finite tuples of predicates. Now we are ready to state the analogue of MacDowell--Specker Theorem which works in the context of the local collection.

\begin{theorem} \label{th_models_with_collection_have_semiregular_end_extensions}
	Let $M$ be a countable piecewise recursively saturated model over a countable signature which satisfies $\PA$ and the full collection scheme. Then there exists a proper elementary end-extension $M \prec_e M'$ such that $M$ is a locally semiregular cut of $M'$.  
\end{theorem}

In the proof, we will use a simple observation which is very far from new or original. It states the equivalence between the collection scheme and the so-called regularity scheme.

\begin{lemma}[Compression lemma] \label{lem_compression}
	Let $M$ be a model in a language with a linear ordering $\leq$ without a largest element, which satisfies the full collection scheme. Let $\phi$ be a binary formula, let $a \in M$, and suppose that there exist arbitrarily large $y \in M$ such that for some $x<a$, $\phi(x,y)$ holds. Then there exists $x_0 <a$ such that $\phi(x_0,y)$ holds for arbitrarily large $y \in M$. 
\end{lemma}
\begin{proof}
	Using the same notation as in the statement of the lemma, assume that for an arbitrarily large elements $y$, there exists an $x<a$ with $\phi(x,y)$. Suppose that for any $x<a$ there exists $z \in M$ such that for any $y>z$, $\neg \phi(x,z)$ holds. Let:
	\begin{displaymath}
		\psi(x,z) \equiv \forall y \Big(y>z \rightarrow \neg \phi(x,y)\Big). 
	\end{displaymath}
	By assumption, for any $x< a$, there exists $z \in M$ such that $M \models \psi(x,z)$. By collection, there exists $b$ such that for any $x<a$, there exists $z<b$ for which $M \models \phi(x,z)$. In particular for any $x<a$ and any $y>b$, 
	\begin{displaymath}
		M \models \neg \phi(x,y),
	\end{displaymath}
contrary to our assumptions. 
\end{proof}

Now we are ready to prove the main result of this section.

\begin{proof}[Proof of Theorem \ref{th_models_with_collection_have_semiregular_end_extensions}]
	The proof is an elaboration of the omitting types argument. Let $M$ be a countable model satisfying full collection and let $a_0,a_1, \ldots$ be an enumeration of its elements. Let $c$ be a fresh constant.	We fix an enumeration of all sentences $\phi_0, \phi_1, \ldots$ in the language expanded with constants $c_0,c_1, \ldots$ and $d_0,d_1,\ldots$. We assume that the constants $c_i, d_i$ do not appear in the sentences $\phi_j$, for $j <i$. The constants $c_i$ are intended as Henkin constants. The constants $d_i$ will play a slightly different role as bounds which allow us to satisfy the local semiregularity condition.
	
	We construct a chain of theories $T_i, i < \omega$. At each step we will construct  auxiliary theories $T_i^1,T_i^2, T_i^3, T_i^4$. Throughout the construction, we assume that all the theories are finite except for the addition of finitely many schemes of the form:
	\begin{displaymath}
		d > a: a \in M.
	\end{displaymath} 
	In effect, we can state of finitely many elements that they are above $M$. 
	
	Let $T_0$ be $\ElDiag(T)$ together with all the sentences of the form $c > a$, $a \in M$. Now suppose that we have defined a theory $T_i$. 
	
	\paragraph*{Step 1}	In order to construct $T^1_{i}$, consider the sentence $\phi_i$. If $T_{i} + \phi_i$ is consistent, we let $T^1_{i}$ be that theory. If not, let $T^1_i = T_i + \neg \phi_i$.
	
	\paragraph*{Step 2}	In order to construct $T^2_i$, consider again the sentence $\phi_i$. If it was not added to the theory constructed in the previous step, let $T^2_i$ be equal to $T^1_i$. If it was, and it has the form $\exists v \psi(v)$, we add the sentence $\psi(c_i)$ to our theory. Note that $c_i$ does not appear in $T_i$ by our bookkeeping assumption. 
	
	\paragraph*{Step 3} In the third step, we check whether the set of sentences of the form 
	\begin{displaymath}
		c_i > a, a \in M
	\end{displaymath}
	 is consistent with $T^2_i$. If yes, then we let $T^3_i$ be obtained from $T^2_i$ by adding this set. Otherwise, there exists some $b \in M$ such that
	 \begin{displaymath}
	 	T^2_i \vdash c_i <b.
	 \end{displaymath} 
	We claim that there exists $d \in M$ such that $T^2_i$ is consistent with the sentence $c_i = d$.
	
	Suppose otherwise. By assumption $T^2_i$ is $\ElDiag(M)$ together with finitely many types of the form $c_j>a, a \in M$ together with finitely many additional formulae. By considering the minimum of the elements $c_j$ and coding the whole tuple as a single element, we can assume without loss of generality that $T^2_i$ extends $\ElDiag(M)$ by a single formula $\psi(t,c_i)$ and a single theory $t>a: a \in M$. (The constant $c_i$ could be of course also eliminated, but we keep it for clarity of the rest of the argument.)  By assumption $T^2_i \vdash c_i < b$ for some $b \in M.$
	
	Now, suppose that a theory $t> a: a \in M  + \ElDiag(M) + \psi(t,c_i)$ proves for any $d \in M$ that $c_i \neq d$. Notice, however, that this implies that for any $d < b \in M$, there exists $a \in M$ such that:
	\begin{displaymath}
		M \models t>a \rightarrow \neg \psi(t,d).
	\end{displaymath} 
	Then, by Compression Lemma, there would exists a single element $r \in M$ for which:
	\begin{displaymath}
		M \models \forall x < b \forall y > r \ \neg \psi (x,y).
	\end{displaymath}
	However this contradicts our assumption on $\psi$. Therefore, there exists in $M$ an element $d <b$ such that for any $a$, we can find $t \in M$ with:
	\begin{displaymath}
		M \models  \psi(d,t) \wedge t>a.
	\end{displaymath}
	Let us pick any such $d$ and set $T^3_i : = T^2_i +  c_i = d.$
	
	\paragraph*{Step 4} Finally, we want to ensure local semiregularity. Suppose that $\phi_i$ is a sentence expressing that some $c_k, k \leq i$ codes a function with the domain $[0,a]$ for some $a \in M$. Specifically:
	
\begin{displaymath}
	\phi_i = \forall x \leq a \exists! y \Big(\tuple{x,y} \in c_k \Big).
\end{displaymath}
In the rest of the argument, let us denote $c_k$ with $f$ for cleaner presentation. We will also replace $\tuple{x,y} \in f$ with the more usual notation $f(x) = y$ and treat $f(x)$ as if it were an independent term. 

\paragraph*{Claim} There exists $a' < a$ and $b \in M$ such that the following theory is consistent (which clearly ensures local semiregularity):
\begin{displaymath}
	T^3_i + \forall x <a'  \Big(  \big(f(x) < b \vee  f(x) > d_i\big) \Big) + d_i >p: p \in M.
\end{displaymath}

We begin with a bit of notation. Recall that $T^3_i$ extends $\ElDiag(M)$ with a sentence $\eta(c_0,\ldots,c_i,d_0, \ldots, d_{i-1} )$ and finitely many sets of sentences of the form $c_j, d_j >q, q\in M$. By considering the minimum of all the parameters and via coding of tuples, we may assume that $T^3_i$ extends $\ElDiag(M)$ simply with a theory of the form:
\begin{displaymath}
	\eta(f,t) + t> q: q \in M,
\end{displaymath}
where $\eta$ is a single formula. Notice that by assumption $\eta(f,t)$ implies that $f$ codes a function with the domain $[0,a]$ in the sense explained above. 

Now, let us say that an interval $[0,r]$ \df{eventually visits} an interval $[b,p]$  if the following holds:
\begin{displaymath}
	\exists N \forall q> N \forall t,f \exists x <r \Big( \eta(t,f) \wedge t>q \rightarrow f(x) \in [b,p] \Big).
\end{displaymath}
We denote it with $\xi(r,b,p)$. 

\paragraph*{Subclaim}
There exists a nonstandard $a_0 \in M$ such that for some $B$, $[0,a_0]$ does not eventually visit any interval $[b,p]$ with $b> B$.

To prove the Subclaim, notice first that for any $n \in \omega$, there exists $B$ for which $\xi(n,b,p)$ does not hold for any $B>n$. Indeed, if no such $B$ exists, then by successively taking intervals which are located higher up in the model, there exists a family of intervals 
\begin{displaymath}
	[b_0,p_0] < [b_1,p_1] < \ldots < [b_n,p_n] < [b_{n+1},p_{n+1}]
\end{displaymath} 
such that $[0,n]$ eventually visits all these intervals. However, this is impossible by the pigeonhole principle, since the image of $[0,n]$ under $f$ cannot intersect $n+2$ disjoint sets. 

Now, we have just shown that the following set of formulae is a type over $M$:
\begin{displaymath}
	\exists B \forall b,p >B \ \neg \xi(v,b,p) \wedge v>\num{n}: n \in \omega.
\end{displaymath}
This type actually uses only finitely many symbols from the signature, so by piecewise recursive saturation, it is realised in $M$. By fixing $a'$ as an element realising this type, we prove the Subclaim. 

We are now ready to prove the claim. Let $a'$ be any element satisfying the subclaim with a bound $b'$. Fix any $p, q \in M$. We claim that there exist $t,f,d_i \in M$ such that:
 \begin{displaymath}
 M \models 	\eta(t,f) \wedge d_i >p \wedge t> q \wedge \forall x <a'  \Big(  \big(f(x) < b \vee  f(x) > d_i\big) \Big)  
 \end{displaymath}
Indeed, by assumption there is no interval $I$ above $b$ such that $[0,a']$ eventually visits $I$. In particular, there exist arbitrarily large $t$ such that $\eta(t,f)$ holds and 
\begin{displaymath}
	\forall x <a' \ f(x) \notin [b,p+1].
\end{displaymath}
 By fixing any such $f$ and $t>q$, and setting $d_i = p+1$, we find our desired interpretation of these three constants in $M$ such that our fixed finite portion of the theory is satisfied. Let $T_i^4 (=: T_i)$, be the theory:
\begin{displaymath}
	T^3_i + \forall x <a'  \Big(  \big(f(x) < b' \vee  f(x) > d_i\big) \Big) + d_i >p: p \in M,
\end{displaymath}
with $a',b'$ chosen as in the proof. (If $\phi_i$ was not of the specific form considered, of course we set $T^4_i = T^3_i$.)

The rest of the argument is completely routine. Let $T_{\omega} = \bigcup_{i \in \omega} T_i$. Consider the Henkin model $N$ given by $T_{\omega}$. In the step for the theories $T^3_i$, we have ensured that every element of $N$ is either greater than all the elements of $M$ or one of these elements. In the fourth step, we have ensured that any function on a domain bounded in $M$ which is definable in $N$ can be restricted to one which has all values either smaller than a fixed element of $M$ or greater than a fixed element of $N$ which ensures local semiregularity. 
\end{proof}
	
	In the above prove, we have used the specific facts that our models admit some degree of saturation which also involved some coding of function. Admittedly, all these assumptions seem somewhat tangential to the main idea of the stated result. In fact, one could imagine another property which makes sense in a more general context. 	Let $M$ be any model over the signature with a binary symbol $<$ interpreted as a linear order. Let $I$ be an initial segment of $M$. We say that $I$ has a \df{gap property} in $M$ if for any formula $\phi(x,y)$ and any $a \in I$, if
	\begin{displaymath}
		M \models \forall x < a \exists y  \ \phi(x,y),
	\end{displaymath} 
	then there exist $b \in I, c \in M$ such that
	\begin{displaymath}
		M \models \forall x < a \exists y \ \phi(x,y) \wedge (y<b \vee y > c).
	\end{displaymath}
It is an interesting question whether any countable model in a signature with a linear order satisfying full collection has an elementary end-extension with a gap property. However, we were unable to obtain our result in this generality.

\section{Questions}

There are some natural questions which we were trying to settle in this article, but which we had to leave open. Let us sum them up.

\begin{enumerate}
	\item We have shown that any countable model of $\CT^-$ with internal collection or internal induction end-extends to a model of the same theory. Our argument crucially uses countability of the model. It seems that a more general proof would need to be much more specifically tailored to the case of truth predicates, since MacDowell--Specker Theorem fails to work in this greater generality.
	
	\emph{Let $M \models \CT^- + \INT (\INTColl)$ be an arbitrary model. Does it have an end-extension satisfying $\CT^- + \INT (\INTColl)$?}
	
	\item With the end-extension results limited to countable models, we are only able to produce $\omega_1$-like models of this theory, in contrast to the general case of models of $\PA$. It seems highly doubtful that $\omega_1$ has indeed such a distinguished status in this context.
	
	\emph{Let $\kappa$ be an arbitrary infinite cardinal. Let $M \models \PA$ be a countable model (alternatively: a model of cardinality $<\kappa$). Does there exist a $\kappa$-like model $M' \succ M$ which expands to a model of $\CT^-$?}\footnote{Ali Enayat has pointed out that a partial answer is provided by Theorem 1.5 in \cite{enayatmohsenipour}, based on  classical results of Chang and Jensen. If $\kappa$ is a cardinal such that $\kappa$ where $\kappa^{<\kappa} = \kappa$ or $\Box_{\kappa}$ holds and $U$ is a consistent theory in a countable language proving full collection, then $U$ has a $\kappa^+$ model. By taking the theory $\ElDiag(M) + \CT^- + \Coll$, we see that under certain set-theoretic assumptions for any countable $M \models \PA$, some $\lambda$-like elementary extensions will exist for $\lambda > \omega_1$.}
	
	\item Our main proof in section \ref{sec_internal_collection} shows that if $M$ is a countable model of a theory in a countable language which extends $\PA$ and features full collection, then it can be end-extended in such a way that it is locally semiregular in the larger model. However, it seems unlikely that this property really depends on anything else but collection (or that wee really have to replace semiregularity with local semiregularity). In particular, the use of coding of sets seems to be a technical artefact of our argument. When coding is removed, semiregularity or its local version no  longer make sense and one has to reword the whole statement which we propose to do in terms of the gap property. It seems likely that collection itself is sufficient to provide end-extensions in which the original model satisfies this modified claim. 
	
	\emph{Let $M$ be a countable model in a signature with a linear order satisfying full collection. Does there exist an elementary end-extension $M' \succ_e M$ such that $M$ satisfies the gap property in $M'$? } 
\end{enumerate}

\section*{Appendix -- a glossary of technical notions}

In the article, we used a number of formalised syntactic notions. For the convenience of the reader, let us gather them in a single glossary. 

\begin{itemize}

	\item $\ClTerm_{\LPA}(x)$ is a formula naturally representing the set of closed terms of $\LPA$.
	\item $\ClTermSeq_{\LPA}(x)$ is a formula naturally representing the set of coded sequences of closed terms of $\LPA$. 
	\item $\dpt(x)$ is a formula naturally representing the depth of a formula $x$. It is a binary relational formula which is provably functional and thus written using the fucntional notation, in accordance with our conventions.
	\item $\form_{\LPA}(x)$ is a formula naturally representing the set of arithmetical formulae. By $\form^{\leq 1}(x)$ we mean the set of formulae with at most one variable free. 		\item $\Sent_{\LPA}(x)$ is a formula naturally representing the set of arithmetical sentences.
	\item $\val{t}$. If $t$ is (a code of) an arithmetical term and $\alpha$ is a $t$-assignment, then by $t^{\alpha}$ or $t(\alpha)$ we mean the formally computed value of $t$ under $\alpha$. If $\ClTerm_{\LPA}(t)$ holds and $\alpha$ is an empty assignment, we write $\val{t}$ instead.  	
\end{itemize}

\section*{Acknowledgements}
We are very grateful to Ali Enayat and Roman Kossak for a number of illuminating discussions, especially concerning classic results in models of $\PA$. This article builds, among other things, on ideas presented in \cite{WcisloKossak} where the crucial role of classic conditions on cuts for the theory of end-extensions of truth predicates first became apparent to us. As we have already mentioned, the article crucially builds upon the ideas worked out by Fedor Pakhomov in his unpublished note \cite{pakhomov_copying_satisfaction}.

This research was supported by an NCN MAESTRO grant 2019/34/A/HS1/00399 ``Epistemic and Semantic Commitments of Foundational Theories.''

\end{document}